\documentclass[12pt]{amsart}
\usepackage{amssymb}

\theoremstyle{plain}
\newtheorem{theorem}{Theorem}
\newtheorem*{lemma*}{Lemma}
\theoremstyle{remark}
\newtheorem*{remark*}{Remark}

\renewcommand{\le}{\leqslant}
\newcommand{\gm}{\Gamma(E)}
\DeclareMathOperator{\Spec}{Spec}
\DeclareMathOperator{\supp}{supp}

\begin{document}
\title[Dual complex]{A note on the dual complex associated to a 
resolution of singularities}
\author{D.~A. Stepanov}
\address{The Department of Mathematical Modelling \\
Bauman Moscow State Technical University \\
Moscow 105005, Russia}
\email{dstepanov@bmstu.ru}
\thanks{The research was supported by RFBR grant no. 05-01-00353}
\date{}

\maketitle

\section{Construction}
The dual complex associated to a resolution of singularities
generalizes the notion of a resolution graph of a surface 
singularity to any dimension. Let $X$ be an algebraic variety over
a field of characteristic $0$ or an analytic space and let 
$S\subset X$ be the singular locus of $X$. Usually we consider
the germ of singularity $(X,S)$. Take a resolution 
$\pi\colon (Y,E)\to (X,S)$ such that the exceptional set $E$ is a
divisor with simple normal crossings. This resolution exists by
Hironaka theorem (see \cite{Hironaka}, \cite{Wlodarczyk}); 
following terminology of the minimal model theory, we say that
$\pi$ is \emph{a log-resolution}.

\emph{The dual complex} $\gm$ \emph{associated to the resolution} 
$\pi$ is a dual $CW$-complex associated to the divisor $E$ as to
a reducible veriety. More precisely, decompose the devisor $E$ to
its prime components $E_i$: $E=\sum_{i=1}^{N}E_i$. Then 
$k$-dimensional cells $\Delta_{i_0\dots i_k}^{j}$ of the complex 
$\gm$ are in one-to-one correspondence with irreducible components 
$E_{i_0\dots i_k}^{j}$ of the intersections
$$E_{i_0}\cap\dots\cap E_{i_k}=\cup_j E_{i_0\dots i_k}^{j}\,,$$
$$k=0,\dots,n-1,\quad 1\le i_0<i_1<\dots<i_k\le N\,.$$
Every cell $\Delta_{i_0\dots i_k}^{j}$ is a standard $k$-dimensional
simplex with vertices marked by numbers $i_0,\dots,i_k$. 
For instance, the $0$-dimensional skeleton of the complex $\gm$ is
a union of the points $\Delta_i$ corresponding to the prime
divisors $E_i$. The patching map from the boundary
$\partial\Delta_{i_0\dots i_k}^{j}$ of a $k$-dimensional cell
to the $k-1$-dimensional skeleton is simplicial. It maps the face
$\Delta_{i_0\dots\widehat{i_s}\dots i_k}^{j}$ onto simplex
$\Delta_{i_0\dots\widehat{i_s}\dots i_k}^{j'}$ such that
$$E_{i_0\dots i_k}^{j}\cap E_{i_0\dots\widehat{i_s}\dots i_k}^{j'}
\ne\varnothing\,,$$
and sends the vertex $i_r$ of the face 
$E_{i_0\dots\widehat{i_s}\dots i_k}^{j}$ to the vertex $i_r$ of the 
simplex $\Delta_{i_0\dots\widehat{i_s}\dots i_k}^{j'}$.

The complex $\gm$ was considered by G.~L. Gordon in the paper
\cite{Gordon}.

\section{Homotopy type of the dual complex is an invariant of
a singularity}
\begin{theorem}\label{T:invariant}
Let $\pi'\colon (Y',E')\to (X,o)$ and $\pi''\colon 
(Y'',E'')\to (X,o)$
be two log-resolutions of an isolated singularity $(X,o)$. Then the 
topological spaces $\Gamma(E')$ and $\Gamma(E'')$ have the same 
homotopy type. 
\end{theorem}
This result is based on the following theorem by 
Abramovich-Karu-Matsuki-W{\l}odarczyk on factorization of
birational maps (see \cite{Matsuki}). 
\begin{theorem}[Weak Factorization Theorem in the Logarithmic
Category]\label{T:factorization}
Let $(U_{X_1},X_1)$ and $(U_{X_2},X_2)$ be complete nonsingular
toroidal embeddings (this means that $U_{X_i}$ is an open set in
a nonsingular variety $X_i$ and the boundary $X_i\setminus U_{X_i}$ 
is a divisor with simple normal crossings, $i=1,2$) over an
algebraically closed field of characteristic zero. 
Let 
$$\varphi\colon(U_{X_1},X_1)--\to (U_{X_2},X_2)$$
be a birational
map which is an isomorphism over $U_{X_1}=U_{X_2}$. Then the map
$\varphi$ can be factored into a sequence of blowups and blowdowns
with smooth admissible and irreducible centers disjoint from
$U_{X_1}=U_{X_2}$. That is to say, there exists a sequence of
birational maps between complete nonsingular toroidal embeddings
$$(U_{X_1},X_1)=(U_{V_1},V_1)\overset{\psi_1}{-\,-\to}(U_{V_2},V_2)
\overset{\psi_2}{-\,-\to}$$
$$\dots\overset{\psi_{i-1}}{-\,-\to}(U_{V_i},V_i)
\overset{\psi_i}{-\,-\to}
(U_{V_{i+1}},V_{i+1})\overset{\psi_{i+1}}{-\,-\to}\dots$$
$$\overset{\psi_{l-2}}{-\,-\to}(U_{V_{l-1}},V_{l-1})
\overset{\psi_{l-1}}{-\,-\to}(U_{V_l},V_l)=(U_{X_2},X_2)\,,$$ 
where \par
(i) $\varphi=\psi_{l-1}\circ\psi_{l-2}\circ\dots\circ\psi_1$, \par
(ii) $\psi_i$ are isomorphisms over $U_{V_i}$, and \par
(iii) either $\psi_i$ or $\psi_{i}^{-1}$ is a morphism obtained
by blowing up a smooth irreducible center $C_i$ (or $C_{i+1}$) 
disjoint from $U_{V_i}=U_{V_{i+1}}$ and transversal to the boundary
$D_{V_i}=V_i\setminus U_{V_i}$ (or $D_{V_{i+1}}=
V_{i+1}\setminus U_{V_{i+1}}$), i. e., at each point $p\in V_i$
 (or $p\in V_{i+1}$) there exists a regular coordinate system
$\{x_1,\dots,x_n\}$ in a neighborhood $p\in U_p$ such that 
$$D_{V_i}\cap U_p\; (\text{or } D_{V_{i+1}}\cap U_p)=
\{\prod_{j\in J}x_j=0\}$$
and 
$$C_i\cap U_p\; (\text{or }C_{i+1}\cap U_p)=\{\prod_{j\in J}x_j=0\,,
\; x_{j'}=0\; \forall j'\in J'\}\,,$$
где $J,J'\subseteq\{1,\dots,n\}$.
\end{theorem}
\begin{remark*}
The same holds if $X_i$ are complex manyfolds and $\varphi$, 
$\psi_i$ are bimeromorphic maps.
\end{remark*}

To apply Theorem~\ref{T:factorization} in our case, take the 
resolutions $(Y'\setminus E',Y')$ and $(Y''\setminus E'',Y'')$ as 
toroidal embeddings and compactify $Y'$ and $Y''$ to smooth
verieties (here we use the fact that the given singularity $(X,o)$
is isolated). Now Theorem~\ref{T:invariant} follows from the 

\begin{lemma*}
Let $\sigma\colon(X'\setminus E',X')\to(X\setminus E,X)$ be a
blowup of an admissible center $C\subset E$ in a nonsingular
toroidal embedding $(X\setminus E,X)$, 
$X'\setminus E'\simeq X\setminus E$. Then the topological spaces
$\Gamma(E')$ and $\Gamma(E)$ have the same homotopy type.
\end{lemma*}
\begin{proof}
Let $E=\sum_{i=1}^{N}E_i$ be the decomposition of $E$ into its
prime components, and let $C\subset E_i$ for $1\le i\le l$ and
$C\nsubseteq E_i$ for $l<i\le N$. Assume that $C$ has nonempty
intersections also with $E_{l+1},\dots,E_r$, $l<r\le N$. There are
two possibilities. \newline
1) $\dim C=n-l$ ($n=\dim X$), i. e., $C$ coincides with one of the
irreducible components of the intersection $E_1\cap\dots\cap E_l$:
$C=E_{1\dots l}^{1}$. Then after the blowup the intersection of the 
proper transforms $E'_1,\dots,E'_l$ of the divisors $E_1,\dots,E_l$ 
has $J-1$ irreducible components (if $J$ is the number of components
of $E_1\cap\dots\cap E_l$), but all these proper transforms
intersect the exceptional divisor $F$ of the blowup $\sigma$.
Futhermore, $F$ intersects proper transforms $E'_{l+1},\dots,
E'_r$ of the divisors $E_{l+1},\dots,E_r$. Now it is clear that
the complex $\Gamma(E')$ is obtained from $\gm$ by the barycentric
subdivision of the simplex $\Delta_{1\dots l}^{1}$ with the center
at the point corresponding to the divisor $F$. Thus the complexes
$\Gamma(E')$ and $\gm$ are even homeomorphic. \newline
2) $\dim C<n-l$, let $C\subset E_{1\dots l}^{1}$. In this case
divisors $E_{i_1},\dots,E_{i_s}$ have nonempty intersection if
and only if their proper transforms $E'_{i_1},\dots,E'_{i_s}$ have
nonempty intersection. Therefore the complex $\Gamma(E')$ is
obtained from the complex $\gm$ in the following way. Add to $\gm$ 
a new vertex corresponding to the exceptional divisor $F$ of the
blowup $\sigma$ and construct cones with vertex at $F$ over all
the maximal cells $\Delta_{i_1\dots i_s}^{j}$ of the complex $\gm$
possessing the property
$$E_{i_1\dots i_s}^{j}\cap C\ne \varnothing\,.$$
Note that the simplex $\Delta_{F,1\dots l}$ corresponding to
the intersection $F\cap E_{1\dots l}^{1}$ is regarded as
a common simplex for all constructed cones. Now we can define the
homotopy equivalence between $\Gamma(E')$ and $\gm$ as a contraction
of the constructed cones: it sends the vertex $F$ of the complex 
$\Gamma(E')$ to any of the vertices $E_1,\dots,E_l$ of the cell
$\Delta_{1\dots l}^{1}$ of the complex $\gm$ and it is identity
on other vertices of $\Gamma(E')$ ($\gm$). Then the induced
simplicial map is our homotopy equivalence. 
\end{proof}

\section{Toric singularities}
As an example let us find the homotopy type of the dual complex 
associated to a resolution of a toric singularity.
\begin{theorem}
Let $o\in X$ be an isolated singularity of a toric variety $X$. 
Then for any log-resolution $\pi\colon(Y,E)\to(X,o)$ the dual 
complex $\gm$ is homotopy equivalent to a point.
\end{theorem}
\begin{proof}
We can restrict ourselves to the case when the variety $X$ is
affine. Then it can be represented as a toric variety
$X=U_{\sigma}=\Spec\mathbb{C}[\sigma]$ corresponding to a convex
polyhedral cone $\sigma\subset\mathbb{R}^n$. According to
Theorem~\ref{T:invariant}, it is enough to determine the homotopy
type for one arbitrary log-resolution of the singularity
$(X,o)$. 

There exists a fan $\Sigma$ such that the support
$\supp\Sigma=\sigma$ and the birational morphism $\pi\colon 
Y=X(\Sigma)\to X$ is a log-resolution for the variety $X$ (see 
\cite{Danilov}). First consider the dual complex for the divisor
$E'=\sum T_i+E$, where $T_i$ are invariant divisors corresponding
to the edges $\tau_i$ of the cone $\sigma$ and $E$ is the 
exceptional divisor of the morphism $\pi$. Recall that the prime
exceptional divisors $E_i$, $E=\sum E_i$, correspond to 
1-dimensional cones of the fan $\Sigma$ different from $\tau_i$; 
the divisors $E_i$ and $E_j$ have a nonempty intersection if and 
only if the corresponding 1-dimensional cones span a 2-dimensional
cone belonging to the fan $\Sigma$, and so on.

Let $L$ be a hyperplane in the space $\mathbb{R}^n$ such that
$L\cap\sigma=K$ is a compact polyhedron. Then the fan $\Sigma$ 
determines some triangulation of the poyhedron $K$. We denote by the
same letter $K$ the corresponding simplicial complex. Now it is
clear that the dual complex $\Gamma(E')$ is homeomorphic to the
complex $K$ and it obviously has the homotopy type of a point.

The dual complex $\gm$ is obtained from $K$ by deleting simplexes
which contain at least one of divisors $T_i$ as a vertex. The
triangulation $\Sigma$ can be chosen sufficiently small, so it
is evident that after such deleting we obtain a complex which is
homotopy equivalent to the original complex $K$, i. e., which has 
the homotopy type of a point.
\end{proof}

\section{Some remarks}
The surface cusp
$$X=\{x^4+y^4+z^4+xyz=0\}\subset\mathbb{C}^3\,;$$
gives an example of singularity with homotopy nontrivial dual
complex associated to a resolution. It is easy to see that if
$\pi\colon Y\to X$ is a log-resolution, $E$ is its exceptional
divisor, then the complex $\gm$ is homotopy equivalent to the
circle $S^1$.

The singularity
$$\{x^8+y^8+z^8+x^2y^2z^2=0\}\subset\mathbb{C}^3$$
gives a more complicated example (see \cite{Gordon}). If we blow up
the origin, the exceptional divisor $E'|_{X'}$ consists of 3 lines
$E_i$, $i=1,2,3$; every 2 of them intersect at a single point. The 
proper transform $X'$ is singular along these lines. A resolution 
of $X$ can be obtained by blowing up the variety $X'$ along 
$E_1\cup E_2\cup E_3$. The dual complex associated to this resolution
has a nontrivial first homology group; in particular,
$\dim H_1(\gm,\mathbb{Q})=4$.

On the other hand, it is known (see \cite{Artin}) that the
resolution graphs of rational surface singularities are trees, thus
they have homotopy type of a point. It would be interesting to
varify if this fact generalizes to higher dimensions. More precisely,
if $(X,o)$ is an isolated rational singularity and
$\pi\colon(Y,E)\to(X,o)$ is a log-resolution, then is it right that
the complex $\gm$ has the homotopy type of a point?

\end{document}